\documentclass[11pt,a4paper]{amsart}

\usepackage{amscd,amsmath,amsfonts,amssymb}

\newcommand{\la}{\lambda}

\def\.{\cdot}

\newcommand{\RR}{\mathbb{R}}

\newtheorem{te}{Theorem}

\begin{document}

\noindent G{\'e}ometrie diff{\'e}rentielle~ /~ Differential Geometry
\bigskip

\title[Eigenvalue estimates for the Dirac operator]{Eigenvalue
  estimates for the Dirac operator and harmonic 1--forms of constant length}
%\title[First eigenvalue of Dirac operator]{The first eigenvalue of the
 % Dirac operator on compact  spin manifolds with non zero first Betti number}
\thanks{The authors are  members of EDGE, Research
Training Network HRPN-CT-2000-00101, supported by the European Human
Potential Programme.}
\author{Andrei Moroianu}
\author{Liviu Ornea}
\address{Centre de Math{\'e}mathiques, Ecole Polytechnique, 91128 Palaiseau Cedex, France}
\email{am@math.polytechnique.fr}
\address{University of Bucharest, Faculty of Mathematics,
14 Academiei str., 70109 Bucharest, Romania}
\email{Liviu.Ornea@imar.ro}


\begin{abstract}\hfill
\par\noindent We prove that on a compact $n$--dimensional  spin
  manifold admitting a non--trivial harmonic
$1$--form of constant length, every eigenvalue $\la$ of the Dirac operator
satisfies the inequality $\la^2 \geq \frac{n-1}{4(n-2)}\inf_M
  \mathop{\mathrm{Scal}}$. In the limiting case the  universal cover
  of the manifold is isometric to $\mathbb{R}\times N$ where $N$ is a
  manifold admitting Killing spinors.

\bigskip\bigskip
\begin{center}
{\Large  \bf  Estimations de valeurs propres pour
l'op{\'e}rateur de Dirac et 1--formes
harmoniques de longueur constante}
\end{center}
\bigskip
\noindent \textsc{R{\'e}sum{\'e}.}~

\noindent Nous d{\'e}montrons que toute valeur propre $\la$ de
l'op{\'e}rateur de Dirac d'une vari{\'e}t{\'e} spinorielle
compacte, de dimension $n$, qui admet une $1$--forme harmonique
non--triviale de
longueur constante
v{\'e}rifie l'in{\'e}galit{\'e}  $\la^2 \geq \frac{n-1}{4(n-2)}\inf_M
  \mathop{\mathrm{Scal}}$. Dans le cas limite  le rev{\^e}tement
universel de la
  vari{\'e}t{\'e} est isom{\'e}trique {\`a} $\RR\times N$ o{\`u} $N$ est une
  vari{\'e}t{\'e} admettant des spineurs de Killing.
\end{abstract}
\maketitle

\section{Introduction}

Let $(M^n,g)$ be a compact spin manifold of (real)
dimension $n$. For positive scalar curvature
$\mathop{\mathrm{Scal}}:=\mathop{\mathrm{Scal}}(M,g)$, every
eigenvalue $\la$ of the Dirac operator $D$  satisfy the well--known
Friedrich inequality \cite{fr}:
\begin{equation}\label{fr}
\la^2 \geq \frac{n}{4(n-1)}\inf_M \mathop{\mathrm{Scal}}.
\end{equation}
There exist manifolds on which the  inequality is indeed sharp: the
limiting case is equivalent to the existence of   a  Killing
spinor, \emph{i.e.} a spinor $\Psi$  satisfying the equation
\begin{equation}
\nabla_X\Psi+\frac{\la}{n}X\cdot\Psi=0.
\end{equation}
On the other hand, Hijazi noticed that a manifold admitting a parallel
$k$--form, $k\neq 0,n$, carries no Killing spinors \cite{hij}.
Consequently, Friedrich's inequality cannot be sharp on manifolds with
geometric structures which support parallel forms. Various authors
have obtained sharp improvements of \eqref{fr} on K{\"a}hler and
quaternionic
K{\"a}hler manifolds (see \cite{k1}, \cite{k2}, \cite{ksw}).

Another recent improvement of Friedrich inequality in this
latter direction was found by Alexandrov, Grantcharov and Ivanov. They
proved in \cite{agi} that the existence of a parallel
$1$--form on $M^n$, $n\ge3$, implies the inequality
\begin{equation}\label{bulg}
\la^2 \geq \frac{n-1}{4(n-2)}\inf_M \mathop{\mathrm{Scal}}.
\end{equation}
The universal covering space of the manifolds appearing in the
limiting case was also described.

In this note we generalize the above result showing that \eqref{bulg}
can be derived only from the existence of a \emph{harmonic} $1$--form
\emph{with constant length}:

\begin{te}\label{main}
Inequality \eqref{bulg} holds on any compact
spin manifold  $(M^n,g)$, $n\ge3$, admitting a non--trivial harmonic
$1$--form $\theta$ of constant
length. The limiting
case is obtained if and only if $\theta$ is parallel and
the eigenspinor $\Psi$ corresponding to the smallest
eigenvalue of the Dirac operator satisfies the following Killing type equation:
\begin{equation}\label{kil}
\nabla_X\Psi+\frac{\la}{n-1}(X\cdot\Psi-\langle X,
\theta\rangle\theta\cdot\Psi)=0,
\end{equation}
after rescaling $\theta$ to unit length.
\end{te}

%We then rely on \cite[Theorem 3.1]{agi} to describe the universal
%cover of $(M,g)$.

%Geometrically, the existence of a harmonic $1$--form of constant length
%is equivalent with the existence of a foliation with minimal
%hypersurfaces whose normal trajectories are geodesics. Moreover,  if
%the form is integral, there
%exists a Riemannian submersion  of the
%manifold over $S^1$.

The condition for the norm of the 1--form to be constant is
essential, in the sense that the topological constraint alone
-- the existence of a non--trivial harmonic 1--form -- does not
allow any improvement of Friedrich's inequality.

Indeed, motivated by a conjecture appearing in an earlier version of this
note, B{\"a}r and Dahl \cite{bd}
have constructed, on any compact spin manifold $M^n$
and for every positive real number $\epsilon$, a metric $g_\epsilon$ on
$M$ with the property that  ${\mathrm{Scal}_{g_\epsilon}}\ge n(n-1)$ and
such that the first eigenvalue of the Dirac operator
satisfies $\lambda_1^2(D_\epsilon)\le \frac{n^2}{4}+\epsilon$. This
construction clearly shows that no improvement of Friedrich's inequality
can be obtained under purely topological restrictions.

{\bf Acknowledgment.} L.O. thanks the Centre of Mathematics of the
Ecole Polytechnique (Pa\-lai\-seau) for hospitality during March--May 2003
when this research was initiated.

\section{The main inequality}

Let $\theta$ be a 1--form of unit length on a spin manifold $(M^n,g)$ and
let $\Psi$ be an arbitrary spinor field on $M$. We identify $1$--forms
with vector fields by means of the scalar product that we denote with
$\langle\  ,\ \rangle$.

Consider the
following ``twistor--like'' operator $T:TM\otimes \Sigma M\to \Sigma M$
$$T_X\Psi=\nabla_X\Psi+\frac{1}{n-1}X\.D\Psi-\frac{1}{n-1}\langle X,
\theta\rangle\theta\.D\Psi-\langle X,\theta\rangle\nabla_\theta\Psi,$$
where the dot $\.$ denotes Clifford multiplication.
A simple calculation yields
\begin{equation}\label{tw}
|T\Psi|^2=|\nabla\Psi|^2-\frac{1}{n-1}|D\Psi|^2-|\nabla_\theta\Psi|^2
+\frac{2}{n-1}\langle D\Psi,\theta\.\nabla_\theta\Psi\rangle.
\end{equation}
From now on we will suppose that $\theta$ is harmonic, $M$ is compact
with volume element $d\mu$ and
has positive scalar curvature $\mathop\mathrm{Scal}$, and $\Psi$ is an eigenspinor of
the Dirac operator $D$ of $M$
corresponding to the least eigenvalue (in absolute value), say $\lambda$.
We let $\{e_i\}$, $i=1,\ldots,n$ denote a  local orthonormal  frame on $M$.

The harmonicity of $\theta$ implies the following useful relation:
\begin{equation}\label{harm}
D(\theta\.\Psi)=-\theta\.D\Psi-2\nabla_\theta\Psi.
\end{equation}
Indeed, one may write:
\begin{equation*}
\begin{split}
D(\theta\.\Psi)&=\sum e_i\.\nabla_{e_i}(\theta\.\Psi)=\sum
e_i\.(\nabla_{e_i}\theta)\.\Psi+e_i\.\theta\.\nabla_{e_i}\Psi\\
&=(d\theta+\delta\theta)\.\Psi + e_i\.\theta\.\nabla_{e_i}\Psi\\
&=-\theta\cdot\sum e_i\cdot\nabla_{e_i}\Psi-2\langle
e_i,\theta\rangle\nabla_{e_i}\Psi\\
&=-\theta\.D\Psi-2\nabla_{\theta}\Psi.
\end{split}
\end{equation*}
Taking the square norm in \eqref{harm} yields
\begin{equation}\label{har}
|D(\theta\.\Psi)|^2=|\theta\.D\Psi|^2+4|\nabla_\theta\Psi|^2
-4\langle D\Psi,\theta\.\nabla_\theta\Psi\rangle,
\end{equation}
%\begin{equation}\label{har}
%\langle D\Psi,\theta\.\nabla_\theta\Psi\rangle=|\nabla_\theta\Psi|^2
%-\frac14(|D(\theta\.\Psi)|^2-|\theta\.D\Psi|^2).
%\end{equation}
By integration over $M$ in \eqref{tw}, using \eqref{har} to express
the last term in the right hand side of \eqref{tw}, and the
Lichnerowicz
formula $D^2=\nabla^*\nabla+\frac 14 \mathop{\mathrm{Scal}}$, we get
\begin{equation}\label{int}
\begin{split}
\int_M|T\Psi|^2d\mu&=\int_M\left\{\frac{n-2}{n-1}|D\Psi|^2-\frac14 \mathop{\mathrm{Scal}}|\Psi|^2-
\frac{n-3}{n-1}|\nabla_\theta\Psi|^2\right.\\
&\left.-\frac{1}{2(n-1)}(|D(\theta\.\Psi)|^2-
|\theta\.D\Psi|^2)\right\}d\mu.
\end{split}
\end{equation}
The term in the last bracket of the integrand is clearly positive since,
from the choice of $\lambda$ to
be minimal, we have from the classical Rayleigh inequality
$$\lambda^2\le \frac{\int_M|D\Phi|^2d\mu}{\int_M|\Phi|^2d\mu}$$
for every $\Phi$. In particular, for $\Phi=\theta\.\Psi$ this reads
$$\int_M|D(\theta\.\Psi)|^2d\mu\ge\lambda^2\int_M|\theta\.\Psi|^2d\mu
=\lambda^2\int_M|\Psi|^2d\mu
=\int_M|D\Psi|^2d\mu=\int_M|\theta\.D\Psi|^2d\mu.$$
Thus \eqref{int} gives
$$\int_M(\frac{n-2}{n-1}\lambda^2-\frac14 \mathop{\mathrm{Scal}})
|\Psi|^2d\mu=
\int_M|T\Psi|^2+\frac{n-3}{n-1}|\nabla_\theta\Psi|^2+
\frac{1}{2(n-1)}(|D(\theta\.\Psi)|^2-|\theta\.D\Psi|^2)
d\mu
\ge0,$$
which immediately implies the first statement of Theorem \ref{main}.

\section{The limiting case}

Suppose now that equality is reached in \eqref{bulg} for the
eigenvalue $\la$ with corresponding eigenspinor $\Psi$. Then
$T\Psi=0$. Contracting with $e_i$~:
$$\sum e_i\.\nabla_{e_i}\Psi+\frac{\la}{n-1}\sum
e_i\.e_i\.\Psi-\frac{\la}{n-1}\sum e_i\.\langle
e_i,\theta\rangle\theta\.\Psi-
e_i\.\langle e_i,\theta\rangle\nabla_\theta\Psi=0,$$
 gives $\theta\.\nabla_\theta\Psi=0$,
so $\nabla_\theta\Psi=0$ (for $n>3$ this follows directly from the
vanishing of the integral $\int_M\frac{n-3}{n-1}|\nabla_\theta\Psi|^2d\mu$).
Thus $\Psi$ satisfies the Killing type equation
\begin{equation}\label{kill}
\nabla_X\Psi=aX\.\Psi-a\langle X,\theta\rangle\theta\.\Psi,\quad a=
-\frac{\la}{n-1}.
\end{equation}
In order to show that $\theta$ is parallel,  we first compute the
spin curvature operator $\mathcal{R}_{Y,X}=[\nabla_Y,\nabla_X]-\nabla_{[Y,X]}$
acting on $\Psi$. We have successively~:
\begin{equation*}
\begin{split}
\frac 1a\nabla_Y\nabla_X\Psi&=\nabla_YX\.\Psi+aX\cdot(Y-\langle
Y,\theta\rangle\theta)\.\Psi-\langle\nabla_YX,\theta\rangle\theta\.\Psi\\
&-\langle X,\theta\rangle\nabla_Y\theta\.\Psi-\langle X,\nabla_Y\theta\rangle\theta\.\Psi-a\langle
X,\theta\rangle\theta\.(Y-\langle
Y,\theta\rangle\theta)\.\Psi\\
\frac 1a\mathcal{R}_{Y,X}\.\Psi&=a(X\.Y-Y\.X)\.\Psi
-a(\langle Y,\theta\rangle X\.\theta-\langle X,\theta\rangle Y\.\theta)\.\Psi\\
&+(\langle Y,\nabla_X\theta\rangle-\langle
X,\nabla_Y\theta\rangle)\theta\.\Psi+(\langle Y,\theta\rangle\nabla_X\theta-\langle
X,\theta\rangle\nabla_Y\theta)\.\Psi\\
&+a(\langle Y,\theta\rangle\theta\.X-\langle
X,\theta\rangle\theta\.Y)\.\Psi\\
&=a(X\.Y-Y\.X)\.\Psi+2a(\langle Y,\theta\rangle\theta\. X-\langle
X,\theta\rangle\theta\.Y)\.\Psi\\
&+(\langle Y,\nabla_X\theta\rangle-\langle
X,\nabla_Y\rangle)\theta\.\Psi+(\langle Y,\theta\rangle\nabla_X\theta-\langle
X,\theta\rangle\nabla_Y\theta)\.\Psi.
\end{split}
\end{equation*}
Using again the harmonicity of $\theta$  we easily
derive~:
\begin{equation}\label{ric}
\begin{split}
\frac{1}{2a}Ric(X)\.\Psi&=\frac 1a\sum e_i\.\mathcal{R}_{e_i,X}\Psi=\\
&=2(n-2)a(X-\langle
X,\theta\rangle\theta)\.\Psi+2\theta\.\nabla_X\theta\.\Psi-2\langle
X,\nabla_\theta\theta\rangle\Psi.
\end{split}
\end{equation}
But $\nabla_\theta\theta=0$ because $\theta$ is unitary and
closed. Indeed, for any vector field $Y$~:
\begin{equation*}
%\begin{split}
\langle
Y,\nabla_\theta\theta\rangle=d\theta(\theta,Y)-\langle\theta,Y\rangle=0.
\end{equation*}
Hence, taking $X=\theta$ in \eqref{ric}, we obtain  $Ric(\theta)=0$.
Then, as $\theta$ is harmonic,  the Bochner formula
assures that $\theta$ is parallel.
This
completes the proof of Theorem \ref{main}.

Since the 1--form $\theta$ has to be parallel in the limiting
case, we can apply Theorem 3.1 in \cite{agi} to determine the
universal cover of $M$. In fact, as proved in \emph{loc. cit.},
this  is isometric to a Riemannian product $\mathbb{R}\times N$,
where $N$ is a spin manifold carrying a real Killing spinor, hence
can be described by  B{\"a}r's classification \cite{bar}. Finally,
$M$ turns out to be a suspension of an isometry of $N/\Gamma$ (a
finite quotient of $N$) over the circle.

\end{document}